\documentclass[a4paper,11pt]{article}
\usepackage{amsmath,amscd,amssymb}
\makeatletter
   \let\accent@spacefactor\relax
 \makeatother
\usepackage[french]{babel}
\usepackage{fancyheadings}
\usepackage{graphicx,psfrag}

\def \mod {{\mathfrak{Mod}}}
\def \ab {{\mathfrak{Ab}}}
\def \F {{\cal F}}
\def \coh {{\mathfrak{Coh}}}
\def \fu #1#2{{\bf\uppercase{#1\footnotesize{#2}}}}
\def \th #1#2#3{\noi\fu{#1}{#2} {\bf #3.} ---}
\def \thf #1#2{\noi\textsc{#1 #2.} ---}
\def \avd {anneau de valuation discr\`ete }

\def \p {{\mathbb P}}

\def \fl {\longrightarrow}
\def \ot {\otimes}
\def \otk {\otimes_{\opk}}
\def \ota {\otimes_{\opa}}

\def \vs {\vskip}

\def \L {{\cal L}}

\def \cL {{\check {\cal L}}}
\def \ccE {{\check {\cal E}}^{\bullet}}

\def \LC {L^{\bullet}}
\def \EC {E^{\bullet}}
\def \ECC {{\check E}^{\bullet}}
\def \CL {{\check L}}

\def \oo {{\cal O}}
\def \dm {{\textit{D\'emonstration}}}

\def \Q {{\mathbb{Q}}}

\def \noi {\noindent}
\def \U {{\bf U}}
\def \ex #1#2{\medbreak\noindent{\bf Exemple #1.#2.}\enspace  }
\def \mapdown#1{\Big\downarrow
\rlap{$\vcenter{\hbox{$\scriptstyle#1$}}$}}
\def \mapsdown#1{\downarrow
\rlap{$\vcenter{\hbox{$\scriptstyle#1$}}$}}
\def \mapup#1{\Big\uparrow
\rlap{$\vcenter{\hbox{$\scriptstyle#1$}}$}}

\def \E {{\cal E}}

\def \op {\oo_{X}}
\def \opk {\oo_{X_k}}

\def \opa {\oo_{X_A}}
\def \opan {\oo_{X_{A_n}}}
\def \opax {\oo_{X_A,x}}
\def \e {\underline{\bf Ext}_{\oo_X}}
\def \ev {\underline{\bf Ext}}
\def \ea {\underline{\bf Ext}_{\opa}}

\def \ek {\underline{\bf Ext}_{\opk}}
\def \h {\underline{\bf Hom}}
\def \ha {\underline{\bf Hom}_{\opa}}

\def \h {\underline{\bf Hom}_{\op}}
\def \hv {\underline{\bf Hom}}
\def \hk {\underline{\bf Hom}_{\opk}}
\def \ta {\underline{\bf Tor}^{\opa}}
\def \t {\underline{\bf Tor}}

\def \tk {\underline{\bf Tor}^{\opk}}
\def \M {{\bf M}}
\def \U {{\bf U}}
\def \bI {{\bf I}}
\def \bV {{\bf V}}
\def \Mg {\mathfrak{M}}

\def \H {\widetilde{H}}

\begin{document}

\centerline{\large{\uppercase{\bf{D\'eformations de fibr\'es
vectoriels}}}}
\centerline{\large{\uppercase{\bf{sur les vari\'et\'es de dimension 3}}}}

\vs 0.3 cm

\centerline{\large{Nicolas PERRIN}}


%
%

\vs 1 cm

\centerline{\large{\bf Introduction}}

\vs 0.4 cm

Soit $X$ une vari\'et\'e projective lisse de dimension 3 sur un corps
$k$ alg\'ebriquement clos.
Nous nous int\'eressons dans cet article \`a l'espace de modules
$\M_{X}(0,c_2,0)$ des faisceaux de rang 2, sans torsion, semi-stables et
de classes de chern $(0,c_2,0)$. Soit $\U$ l'ouvert form\'e des
faisceaux localement libres. Dans le cas o\`u $X$ est $\p^3$, cet
ouvert $\U$ contient lui m\^eme l'ouvert $\bI_{c_2}$ des instantons de
degr\'e $c_2$ (form\'e des fibr\'es $E$ stables et v\'erifiant
l'annulation du groupe $H^1E(-2)$). Nous \'etudions les faisceaux $E$
non localement libres qui sont limites de faisceaux localement libres,
c'est \`a dire le bord de $\U$ dans $\M_{X}(0,c_2,0)$.

Cette \'etude a \'et\'e motiv\'ee par les descriptions
de diff\'erentes composantes du bord des instantons. En particuler,
dans \cite{PE1} proposition 15 nous trouvons une obstruction \`a \^etre
limite d'instantons pour une famille de faisceaux de $\M_{\p^3}(0,3,0)$. De
fa\c con plus g\'en\'erale dans $\M_{\p^3}(0,c_2,0)$, la condition
d'Ellingsrud et Str\o mme (cf. \cite{ES} et \cite{HA}) peut \^etre
interpr\'et\'ee comme une obstruction \`a \^etre limite d'instantons
(voir proposition 3.9).

Dans une premi\`ere partie (deux premiers paragraphes) nous donnons
des conditions d'obstruction g\'en\'erales expliquant ces diff\'erents
exemples. Ces r\'esultats permettent de donner une id\'ee du faisceau
g\'en\'eral d'une composante du bord de la vari\'et\'e $\bI_{c_2}$. En
particulier ils sont en accord avec la conjecture de L. Gruson et
G. Trautmann pr\'edisant le bord de $\bI_3$ (voir \cite{PE1} pour plus de
d\'etails).

Dans une seconde partie, nous donnons un exemple sur $\p^3$ qui montre
que les conditions n\'ecessaires de la premi\`ere partie ne sont pas
vides. De plus dans cet exemple elle sont suffisantes.

\vs 0.2 cm

Soit $E$ un faisceau de $\M_{X}(0,c_2,0)$ sans torsion non localement
libre.
Nous donnons des conditions n\'ec\'essaires pour que ce faisceau $E$
soit limite \textit{r\'eduite} (cf. d\'efinition 1.2) de fibr\'es
vectoriels.
Notons $E''$ le bidual de $E$ et $Q$ le conoyau de l'injection de $E$
dans $E''$, nous montrons plus pr\'ecis\'ement que si $E$ a un lieu
singulier de dimension pure 1 (cf. d\'efinition 1.4) alors

\vs 0.2 cm

\th{T}{h\'eor\`eme}{0.1} \textit{Si $E$ est limite
r\'eduite de fibr\'es vectoriels, alors on a un isomorphisme
sym\'etrique $\e^2(Q,\oo_X)\simeq Q$. Soit $C$ le $0^{\rm i\grave{e}me}$
id\'eal de Fitting de $Q$, si $\omega_X\vert_C$ a une racine carr\'ee
$\sqrt{\omega_X\vert_C}$, alors le faisceau
$Q\otimes\sqrt{\omega_X\vert_C}$ est une th\'eta-caract\'eristique sur
$C$.}

\vs 0.2 cm

Plus g\'en\'eralement (sans hypoth\`ese sur le lieu singulier), nous
montrons le

\vs 0.2 cm

\th{T}{h\'eor\`eme}{0.2} \textit{Si le faisceau $E$ est limite r\'eduite de fibr\'es
vectoriels, alors on a la suite exacte
suivante :
\vs -0.3 cm
$$0\fl Q\stackrel{i}{\fl}\e^1(E,\op)\fl\e^3(Q,\op)\fl 0$$}

\vs -0.3 cm

\newpage

Par ailleurs, pour tout faisceau $\F$ d\'efini sur $\op$, il existe une
fl\`eche naturelle de faisceaux
$\h(E'',\F)\stackrel{\psi}{\fl}\e^2(\e^1(E,\op),\F)$. Notons $p$ la
fl\`eche du faisceau  $\h(Q,\F)$

dans $\h(E'',\F)$ d\'efinie par
la surjection de $E''$ dans $Q$. Nous montrons alors le

\vs 0.2 cm

\th{T}{h\'eor\`eme}{0.3} \textit{Si le faisceau $E$ est limite r\'eduite de
fibr\'es vectoriels, alors la compos\'ee $\ev^2(i)\circ\psi\circ p$ de
$\h(Q,\F)$ dans $\e^2(Q,\F)$ est nulle.}

\vs 0.2 cm

Dans le cas des faisceaux $E$ dont le lieu singulier est de dimension pure
0 construits dans \cite{PE1}, cette condition est suffisante. Cependant, si
le faisceau $E$ a un lieu singulier de dimension pure \'egale \`a 
1, nous donnons une nouvelle une condition. Nous
montrons 
que la surjection de $E''$
dans $Q$ (de support $C$) d\'efini un \'el\'ement $\xi\in{\rm
Ext}_{\oo_C}^1(\e^1(Q,\F),\e^2(Q,\F))$. Nous montrons alors le

\vs 0.2 cm

\th{T}{h\'eor\`eme}{0.4} \textit{Si le faisceau $E$ est limite r\'eduite de
fibr\'es vectoriels, alors l'\'el\'ement $\xi$ du groupe d'extensions
${\rm Ext}_{\oo_C}^1(\e^1(Q,\F),\e^2(Q,\F))$ est nul.}

\vs 0.2 cm

Nous montrons ensuite en utilisant la th\'eorie des foncteurs coh\'erents
(\cite{Ha2}) que les conditions fonctorielles des th\'eor\`emes 0.3 et 0.4
sont v\'erifi\'ees si et seulement si elles le sont pour le seul faisceau
$\F=Q$. Plus pr\'ecis\'ement nous montrons la

\vs 0.2 cm

\th{P}{roposition}{0.5} \textit{Notons $\coh(X)$ la cat\'egorie des
  foncteurs coh\'erents de $\mod(X)$.}

\textit{(\i) Sous les hypoth\`eses du th\'eor\`eme
0.3, on a l'identification :} 
$${\rm Hom}_{\coh(X)}(\h(Q,\bullet),\e^2(Q,\bullet))\simeq\e^2(Q,Q)$$

\textit{(\i\i) Sous les hypoth\`eses du th\'eor\`eme 0.4, on a
l'identification :} 
$${\rm Ext}^1_{\coh(X)}(\e^1(Q,\bullet),\e^2(Q,\bullet))\simeq H^1\e^1(Q,Q).$$
\textit{Ce dernier groupe est la valeur prise par le foncteur de gauche en
  $Q$.}

\vs 0.2 cm

Les conditions obtenues sont des conditions pour qu'un faisceau soit
\textit{limite r\'eduite} de fibr\'es vectoriels. Cependant il est
facile de
v\'erifier que, pour les familles construites dans \cite{PE1}, une
d\'eformation \textit{g\'en\'erale} est r\'eduite. Les th\'eor\`emes
0.2 et 0.3 expliquent alors les conditions obtenues \`a la proposition
15 de \cite{PE1}.

\vs 0.2 cm

Dans la seconde partie de cet article nous donnons un exemple qui
montre que les conditions obtenues aux th\'eor\`emes 0.1 et 0.4 ne sont pas
vides et qu'elles
sont dans ce cas suffisantes. Nous d\'ecrivons une famille $\bV$ de
$\M_{\p^3}(0,c_2,0)$ non adh\'erente \`a $\U$ et nous montrons que le
ferm\'e $\bV_0$ v\'erifiant la condition du th\'eor\`eme 0.4 est dans
l'adh\'erence de $\U$ (th\'eor\`eme 3.2). Nous montrons que tous les
faisceaux de $\bV_0$ peuvent \^etre obtenus par d\'eformation
r\'eduite. Si on note $\bV_0^{pair}$ la partie de $\bV_0$ correpondant aux
th\'eta-caract\'eristiques paires, on a :

\vs 0.2 cm

\th{T}{h\'eor\`eme}{0.6} \textit{La condition du th\'eor\`eme 0.4 d\'efini
  dans $\bV$ un ferm\'e strict $\bV_0$. Ce dernier est adh\'erent \`a
  $\U$. Si $c_2\leq 5$ alors $\bV_0^{pair}$ forme une composante
  irr\'eductible du bord de $\bI_{c_2}$.}

\vs 0.2 cm

Nous montrons que toutes les d\'eformations construites dans cette seconde
partie sont r\'e\-dui\-tes. On peut alors formuler la conjecture suivante :

\vs 0.2 cm

\th{C}{onjecture}{0.7} \textit{Une d\'eformation g\'en\'erale de fibr\'es
vectoriels est r\'eduite.}

\vs 0.2 cm

Cette conjecture est vraie pour toutes les familles consid\'er\'ees
dans \cite{PE1} et pour les d\'eforma\-tions que
nous construisons dans la seconde partie de cet article.
Il est \'egalement facile de v\'erifier
qu'elle est vraie pour $\M_{\p^3}(0,2,0)$ ce qui permet d'expliquer
les diff\'erentes composantes du bord de $\bI_2$ d\'ecrites par
M. S. Narashiman et G. Trautmann dans \cite{NT}.

\vs 0.2 cm

Enfin, pour les familles de faisceaux consid\'er\'ees dans la seconde
partie, celles construites dans \cite{PE1} et celles de \cite{NT}, les
conditions pr\'ec\'edentes sont suffisantes.

\vs 0.2 cm

\noi
{\bf Remerciements} : Je tiens \`a remercier Laurent Gruson pour toute
l'aide qu'il m'a apport\'e durant la pr\'eparation de ce travail.

\section{Pr\'eliminaires}

Nous donnons dans ce paragraphe quelques d\'efinitions et les
premi\`eres propri\'et\'es sur les faisceaux limites de fibr\'es
vectoriels dont nous aurons besoin dans la suite.

\subsection{D\'efinitions et premi\`eres remarques}

Nous d\'efinissons les notions de ``limite'' et ``limite r\'eduite''
d'instantons et de ``dimension pure'' du lieu singulier, nous fixons
\'egalement quelques notations.

\vs 0.4 cm

\thf{Notation}{1.1} Soit $A$ un \avd de corps des fractions $K$ et
de corps r\'esiduel $k$. Nous supposerons toujours que $A$ contient
son corps r\'esiduel.
Si $M$ est un $A$-module on note $M_K$ le
$K$ espace vectoriel $M\ot_A K$ et $M_k$ le $k$ espace vectoriel
$M\ot_A k$. Notons $\pi$ une uniformisante de $A$.
On note $A_n$ l'anneau $A/(\pi^n)$.

\vs 0.4 cm

\th{D}{\'efinition}{1.2} Une \textit{d\'eformation de faisceaux} de
$\M_{X}(0,c_2,0)$ d'anneau $A$ est la donn\'ee d'un $\opa$-module $\E$ de
rang 2, plat sur $A$ et tel que $\E_A$ d\'efinisse un point de
$\M_{X}(0,c_2,0)(A)$. Le faisceau $\E_k$ est \textit{la limite} de la
d\'eformation. On dira que c'est une d\'eformation (ou que $\E_k$
est limite) de faisceaux localement libres (resp. d'instantons) si
$\E_K$ d\'efinit un point de $\U(K)$ (resp. $\bI_{c_2}(K)$).

Une d\'eformation $\E$ de fibr\'es vectoriels sera dite
\textit{r\'eduite}  si pour tout $i>0$, les faisceaux
$\ea^i(\E,\opa)$ --- qui sont a priori d\'efinis sur $\opan$ pour un
certain $n$ --- sont d\'efinis sur $\opk$.

\vs 0.4 cm

\thf{Remarque}{1.3} Nous consid\'erons une famille de faisceaux sur
un anneau de valuation discr\`ete $A$. Il n'est pas \'evident a priori
que toute d\'eformation soit de ce type car l'espace des modules n'a
pas de faisceau universel. Cependant, cet espace est un bon quotient
d'une vari\'et\'e sur laquelle il existe un tel faisceau
(cf. \cite{MA}). Il suffit donc de consid\'erer une d\'eformation sur
cette vari\'et\'e, ce qui d\'efinit un faisceau sur un anneau de valuation
discr\`ete $A$. On peut de plus supposer que l'anneau $A$ contient
son corps r\'esiduel $k$ car
on est au dessus d'un corps $k$ alg\'ebriquement clos.

\vs 0.4 cm

\th{D}{\'efinition}{1.4} On dira qu'un faisceau $E$ de rang 2 sur
$X$ et sans torsion a un lieu singulier de dimension pure \'egale
\`a $1$ si les faisceaux $\e^j(E,\op)$ sont nuls pour $j>1$ et si les
faisceaux $\e^k(\e^1(E,\op),\op)$ sont nuls pour $k\not=2$.

\vs 0.4 cm

\th{D}{\'efinition}{1.5} On dira qu'un faisceau $E$ de rang 2 sur
$\p^3$ et sans torsion a un lieu singulier de dimension pure \'egale
\`a $0$ si les faisceaux $\e^j(E,\op)$ sont nuls pour $j>2$ et
si les faisceaux $\e^k(\e^1(E,\op),\op)$ et
$\e^k(\e^2(E,\op),\op)$ sont nuls pour $k\not=3$.

\vs 0.4 cm

\thf{Remarques}{1.6} (\i) On note $E''$ le bidual de $E$ et $Q$ le conoyau de
l'injection canonique de $E$ dans $E''$. Le faisceau $E''$ est
r\'eflexif tandis que le support du faisceau $Q$ est de dimension au
plus 1. Le faisceau $E$ a un lieu singulier de dimension pure 1 ou 0
si on a respectivement l'une des deux conditions suivantes :
\begin{itemize}
\item $E''$ est localement libre et le $0^{i\grave{e}me}$ id\'eal de
Fitting de $Q$ est de dimension pure \'egale \`a 1.
\item le $0^{i\grave{e}me}$ id\'eal de Fitting de $Q$ est de dimension
0.
\end{itemize}

(\i\i) Les faisceaux dont le lieu singulier est de
dimension pure 0 ou 1 forment deux ouverts disjoints de
$\M_{X}(0,c_2,0)\setminus \U$.
Il est ainsi possible que pour un faisceau $E$ assez g\'en\'eral du
bord de $\U$, le lieu singulier de $E$ soit en dimension pure 0 ou
1. C'est ce qu'affirme la conjecture de L. Gruson et G. Trautmann
dans le cas des instantons (cf. \cite{PE1}).

\subsection{Premi\`eres propri\'et\'es}

\th{P}{roposition}{1.7} \textit{Soit $\E$ une d\'eformation de fibr\'es
vectoriels. Le faisceau $\E$ est r\'eflexif de
dimension projective inf\'erieure ou \'egale \`a 2. Si le lieu
singulier de $\E\ota\opk$ est de dimension pure \'egale \`a 1 alors sa
dimension projective est 1.}

\vs 0.2 cm

\dm :
On a ${\rm pd}(\E_x)={\rm prof}(\opax)-{\rm prof}(\E_x)$
(\cite{Ei} Th\'eor\`eme 19.9). Il suffit donc de montrer que si ${\rm
dim}\opax=i\geq 2$, alors ${\rm prof}(\E_x)\geq 2$ pour obtenir le
premier r\'esultat (cf. \cite{Ha} proposition 1.3).

Si $\E_x$ est libre, sa profondeur est \'egale \`a celle de
$\opax$. C'est le cas d\`es que l'id\'eal maximal $\Mg_x$ ne
contient pas une puissance de $\pi$.

Si $\Mg_x$ contient une puissance $\pi^k$ de $\pi$, remarquons que
comme $\E$ est plat sur $A$, l'\'el\'ement $\pi^k$ n'est pas diviseur
de $0$. De plus le faisceau $\E_k$ est sans torsion donc sa profondeur
est au moins 1 et il en est de m\^eme de $\E\ot_A A/(\pi^k)$. Ainsi la
profondeur de $\E_x$ est au moins 2.

Pour prouver la derni\`ere assertion il suffit de montrer que si ${\rm
dim}\opax=4$, alors on a ${\rm prof}(\E_x)\geq 3$. Soit $x$ un tel
point et supposons que le lieu singulier de $\E\ota\opk$ est de
dimension pure \'egale \`a 1. Il suffit de montrer que ${\rm prof
}((\E_k)_x)=2$. Cependant on a la suite exacte $0\fl (\E_k)_x\fl
(\E''_k)_x\fl Q_x\fl 0$ o\`u ${\rm prof }((\E_k'')_x)=3$ et ${\rm prof
}(Q_x)=1$ (cf. remarque 1.6). Ceci prouve le r\'esultat (\cite{Ei}
corollaire 18.6).

\vs 0.4 cm

\thf{Notations}{1.8} (\i) Tous les complexes consid\'er\'es
sont \`a degr\'es d\'ecroissants et nul \`a droite.

(\i\i) Pour toute d\'eformation $\E$ de fibr\'es
vectoriels, nous fixons les notations suivantes. La proposition
pr\'ec\'edente permet de consid\'erer le complexe
$$\E^{\bullet}\ :\ \L_2\fl\L_1\fl \L_0$$
de faisceaux localements libres sur $\opa$ qui est une r\'esolution de
$\E$. Le complexe dual est :
$${\check \E}^{\bullet}\ :\ \cL_0\fl\cL_1\fl\cL_2$$
dont la cohomologie est $\E$, $\ea^1(\E,\opa)$
et $\ea^2(\E,\opa)$
en degr\'es 2, 1 et 0 respectivement. On consid\`ere
\'egalement le complexe
$$\L^{\bullet}\ :\ \L_2\fl\L_1\fl\L_0\fl\cL_0\fl\cL_1\fl\cL_2$$
de cohomologie $\ea^2(\E,\opa)$
et $\ea^1(\E,\opa)$ en degr\'es 0 et 1 et nulle ailleurs.

(\i\i) On note $E$ le faisceau $\E\ota\opk$ et on consid\`ere les
restrictions des complexes pr\'ec\'edents \`a $\opk$. On les note
respectivement de la fa\c con suivante : le complexe
$$E^{\bullet}\ :\ L_2\fl L_1\fl L_0$$
qui est une r\'esolution localement libre de $E$ sur $\opk$, le complexe
$${\check E}^{\bullet}\ :\ \CL_0\fl\CL_1\fl\CL_2$$
dont la cohomologie est $E''$ (le bidual de $E$), $\ek^1(E,\opk)$ et
$\ek^2(E,\opk)$ en degr\'es 2, 1 et 0 et le complexe
$$L^{\bullet}\ :\ L_2\fl L_1\fl L_0\fl\CL_0\fl\CL_1\fl\CL_2$$
dont la cohomologie est nulle en degr\'es 3, 4 et 5 et vaut
$\ek^2(E,\opk)$, $\ek^1(E,\opk)$ et $Q$ (le conoyau de l'injection
$E\fl E''$) en degr\'es 0, 1 et 2 respectivement.

\vs 0.4 cm

\thf{Remarque}{1.9} Soit $\L$ un $\opan$-module. Si $\L$ est d\'efini
sur $\opk$, alors on a l'identification
$\ta_1(\L,\opk)\simeq\L\ota\opk$.

\vs 0.4 cm

\th{L}{emme}{1.10} \textit{Soit $\L$ un $\opk$-module. Pour tout $i$
positif on a la d\'ecomposition suivante (par convention le
bifoncteur $\e^{-1}(\ \bullet\ ,\ \bullet\ )$ est nul) :}
$$\ea^i(\L,\F)=\ek^i(\L,\F)\oplus\ek^{i-1}(\L,\F)$$

\vs 0.2 cm

\dm :
Soit ${L'}^{\bullet}$ une r\'esolution de $\L$ sur $\opk$. Comme $A$
contient son corps r\'esiduel $k$, on peut consid\'erer le
complexe ${L'}^{\bullet}\ot_k A$ (ce qui en fait un $\opa$-module) et
le mapping c\^one ${\L'}^{\bullet}$ du morphisme de complexes
${L'}^{\bullet}\ot_k A\stackrel{\pi}{\fl}{L'}^{\bullet}\ot_k A$.

\vs 0.2 cm

\th{F}{ait}{\!} \textit{Le complexe ${\L'}^{\bullet}$ est une r\'esolution
de $\L$ sur $\opa$.}

\vs 0.2 cm

\dm :
On a la suite exacte de complexes
$$0\fl {L'}^{\bullet}\ot_kA\fl {\L'}^{\bullet}\fl
{L'}^{\bullet}\ot_kA[-1]\fl 0$$
La suite exacte longue de cohomologie permet de v\'erifier que
$H^3({\L'}^{\bullet})$ et $H^2({\L'}^{\bullet})$ sont nuls et que l'on a une
suite exacte :
$$0\fl H^1({\L'}^{\bullet})\fl \L\ot_kA\stackrel{1\ot\pi}{\fl}\L\ot_kA\fl
H^0({\L'}^{\bullet})$$
Cependant la fl\`eche centrale est injective et son conoyau est
$\L\ot_kk\simeq\L$.

\vs 0.2 cm

La suite exacte de complexe pr\'ec\'edente
n'est pas scind\'ee mais elle se scinde sur $\opk$. On applique le
foncteur $\ha(\ \bullet\ ,\F)$ \`a cette suite exacte de
complexe. Elle reste exacte et se scinde. On en d\'eduit
l'identification :
$$H^i(\ha({\L'}^{\bullet},\F))\simeq H^i(\ha({L'}^{\bullet}\ot_kA,\F))\oplus
H^{i-1}(\ha({L'}^{\bullet}\ot_kA,\F)).$$
Le groupe $H^i(\ha({\L'}^{\bullet},\F))$ est $\ea^i(\L,\F)$. De
plus le complexe ${L'}^{\bullet}$ est une r\'esolution de $\L$ sur
$\opk$ donc le groupe
$H^i(\ha({L'}^{\bullet}\ot_kA,\F))=H^i(\hk({L'}^{\bullet},\F))$ et
le groupe $H^{i-1}(\ha({L'}^{\bullet}\ot_kA,\F))=H^{i-1}(\hk({L'}^{\bullet},\F))$
sont $\ek^i(\L,\F)$ et $\ek^{i-1}(\L,\F)$.

\section{Conditions \`a la limite}

Nous supposons dans ce paragraphe que $E$ est un point de
$\M_{X}(0,c_2,0)$. S'il est limite de fibr\'es vectoriels, nous
noterons $\E$ la d\'eformation associ\'ee, le faisceau $E$
s'identifiant \`a $\E\ota\opk$.

\subsection{Propri\'et\'es de dualit\'e}

\th{L}{emme}{2.1} \textit{Si $E$ est limite de fibr\'es vectoriels, les
faisceaux $Q$ et $\ek^2(E,\opk)$ s'identifient aux faisceaux
$\ta_1(\ea^1(\E,\opa),\opk)$ et $\ea^2(\E,\opa)\otk\opk$. On a
\'egale\-ment la suite exacte : \\
\noi
$0\to\ea^1(\E,\opa)\ota\opk\to\ek^1(E,\opk)\to\ta_1(\ea^2(\E,\opa),\opk)\to
0$}

\vs 0.4 cm

\dm :
On a le complexe de faisceaux localement libres (cf. notations 1.8) :
$$\L^{\bullet}\ :\ \L_2\fl\L_1\fl\L_0\fl\cL_0\fl\cL_1\fl\cL_2$$
dont la cohomologie est $\ea^1(\E,\opa)$ et $\ea^2(\E,\opa)$.

Si on applique le foncteur $\ \bullet\ \ota\opk$ \`a ce complexe,
nous obtenons le complexe
$$\LC\ :\ L_2\fl L_1\fl L_0\fl \CL_0\fl\CL_1\fl\CL_2$$
dont la cohomologie est $Q$, $\ek^1(E,\opk)$ et $\ek^2(E,\opk)$.
On a ainsi une suite spectrale
$$\ ^2E^{p,q}=\ta_p(H^q(\L^{\bullet}),\opk)\Rightarrow H^{p+q}(\LC)$$
qui donne les r\'esultats voulus.

\vs 0.4 cm

\thf{Remarque}{2.2} Si le lieu singulier de $E$ est de dimension
pure 1, on a alors les simplifications suivantes : le faisceau
$\ek^2(E,\opk)$ est nul et les faisceaux $Q$ et $\ek^2(Q,\opk)$
s'identifient \`a $\ta_1(\ea^1(\E,\opa),\opk)$ et \`a
$\ea^1(\E,\opa)\ota\opk$.

\vs 0.4 cm

\th{T}{h\'eor\`eme}{2.3} \textit{Si $E$ est limite
r\'eduite de fibr\'es vectoriels et que son lieu singulier est de
dimension pure 1, alors on a un isomorphisme
sym\'etrique $\e^2(Q,\oo_X)\simeq Q$. Soit $C$ le $0^{\rm i\grave{e}me}$
id\'eal de Fitting de $Q$, si $\omega_X\vert_C$ a une racine carr\'ee
$\sqrt{\omega_X\vert_C}$, alors le faisceau
$Q\otimes\sqrt{\omega_X\vert_C}$ est une th\'eta-caract\'eristique sur
$C$.}

\vs 0.2 cm

\dm :
Le lemme 2.1 et la remarque 1.9 (appliqu\'ee \`a $\L=\ea^1(\E,\opa)$)
donnent l'isomorphisme sym\'etrique $\ek^2(Q,\opk)\simeq Q$.
Si $\sqrt{\omega_X\vert_C}$ existe, 
la suite spectrale de changement de base
donne un isomorphisme sym\'etrique
$\hv_{\oo_C}(Q\ot\sqrt{\omega_X\vert_C},\omega_C)\simeq Q\ot\sqrt{\omega_X\vert_C}$.

\vs 0.4 cm

\th{T}{h\'eor\`eme}{2.4} \textit{Si $E$ est limite r\'eduite de
fibr\'es vectoriels, alors on a la suite exacte
suivante :
$$0\fl Q\stackrel{i}{\fl}\ek^1(E,\opk)\fl\ek^3(Q,\opk)\fl 0$$}

\vs -0.2 cm

\dm :
Ce th\'eor\`eme d\'ecoule directement de la remarque 1.9 (appliqu\'ee
aux fais\-ceaux
$\ea^1(\E,\opa)$ et $\ea^2(\E,\opa)$), du lemme 2.1 et de
l'identification entre le faisceau $\ek^3(Q,\opk)$ et $\ek^2(E,\opk)$. Il montre en
particulier que le lieu (r\'eduit) o\`u le faisceau $E''$ n'est pas
localement libre est exactement situ\'e au niveau du support
(r\'eduit) de $Q$. Ceci reste vrai si la d\'eformation n'est pas
r\'eduite.

\subsection{Propri\'et\'es d'annulation}

Soit $E$ un faisceau de $\M_{X}(0,c_2,0)$ et
supposons que $E$ v\'erifie les
conditions du th\'eor\`eme 2.4. Nous cherchons dans ce
paragraphe \`a d\'eterminer des conditions n\'ecessaires sur la
fl\`eche surjective $p$ de $\p({\rm Hom}(E'',Q))$ d\'efinissant $E$
\`a partir de $E''$ et de $Q$ pour que le faisceau $E$ soit limite de
fibr\'es vectoriels.

\vs 0.2 cm

Appliquons le foncteur $\hk(\ \bullet\ ,\F)$ au complexe $\ECC$
(cf. notations 1.8),
nous obtenons le complexe $\EC\otk\F$ et donc la suite spectrale :
\vs -0.3 cm
$$\ ^2E^{p,q}=\ek^p(H^q(\EC),\F)\Rightarrow \tk_{p+q}(E,\F)\ \ \ \ (S)$$
qui donne la fl\`eche suivante :
$$\hk(E'',\F)\stackrel{\psi}\fl\ek^2(\ek^1(E,\opk),\F)$$

\vs 0.4 cm

\thf{Remarques}{2.5} (\i) Si $E$ est limite (r\'eduite ou non) de
fibr\'es vectoriels,
il existe une fl\`eche  $Q\stackrel{i}{\fl}\ek^1(E,\opk)$ qui est la
compos\'ee du morphisme $Q\fl\ea^1(\E,\opa)\ota\opk$ (c'est un
isomorphisme si $E$ est limite r\'eduite) avec la fl\`eche de
restriction canonique

\noi
$\ea^1(\E,\opa)\ota\opk\fl\ek^1(E,\opk)$.

\vs 0.3 cm

(\i\i) Soit ${\cal G}$ un $\opa$-module et $\F$ un $\opk$-module, la
suite spectrale de changement de base
\vs -0.2 cm
$$\ek^p(\ta_q({\cal G},\opk),\F)\Rightarrow\ea({\cal G},\F)$$
donne des fl\`eches fonctorielles \textit{injectives} suivantes :
$$g_p({\cal G})\ :\ \ek^p({\cal
G}\ota\opk,\F)\hookrightarrow\ea^p({\cal G},\F)$$

(\i\i\i) On a les fl\`eches fonctorielles suivantes :
$$f_{p,q}\ :\ \ek^p(\ek^q(E,\opk),\F)\fl\ea^p(\ea^q(\E,\opa),\F)$$
En effet, comme $\E$ est plat on a les fl\`eches de restriction
$\ea^q(\E,\opa)\to\ek^q(E,\opk)$ qui donnent
$\ea^p(\ek^q(E,\opk),\F)\stackrel{{\rm
res}_A}{\fl}\ea^p(\ea^q(\E,\opa),\F)$. Par ailleurs,
la remarque pr\'ec\'edente ou le lemme 1.10 (avec pour $\L$ le
faisceau $\ek^q(E,\opk)$) donnent les fl\`eches :
$\ek^p(\ek^q(E,\opk),\F)\stackrel{g_p}{\fl}
\ea^p(\ek^q(E,\opk),\F)$.

\vs 0.4 cm

\th{T}{h\'eor\`eme}{2.6} \textit{Si $E$ est limite r\'eduite de fibr\'es
  vectoriels, alors on a l'annulation de la compos\'ee}
$$\hk(Q,\F)\stackrel{p}{\fl}\hk(E'',\F)\stackrel{\psi}{\fl}
\ek^2(\ek^1(E,\opk),\F)\stackrel{i}{\fl}\ek^2(Q,\F)$$

\vs 0.2 cm

\dm :
En appliquant le foncteur $\ha(\
\bullet\ ,\F)$ au complexe ${\check \E}^{\bullet}$
on obtient :
\vs -0.3 cm
$$\ ^2E^{p,q}=\ea^p(H^q({\check \E}^{\bullet}),\F)\Rightarrow
\ta_{p+q}(\E,\F).$$
Elle donne les fl\`eches compatibles suivantes :

$$\begin{array}{ccc}
\ha(\ha(\E,\opa),\F) & \fl & \ea^2(\ea^1(\E,\opa),\F) \\
\uparrow & & \mapup{f_{2,1}} \\
\hk(\hk(E,\opk),\F) & \fl & \ek^2(\ek^1(E,\opk),\F) \\
\end{array}$$

\vs 0.2 cm

Le faisceau $\hk(E,\opk)$ est isomorphe \`a $E''$ et
$\ha(\ha(\E,\opa),\F)$ est isomorphe \`a $\hk(E,\F)$ (car $\E$ est
autodual, plat sur $A$ et que $\F$ est d\'efini sur $\opk$).

La suite
exacte $\hk(Q,\F)\fl\hk(E'',\F)\fl\hk(E,\F)$ impose que l'image du
faisceau $\hk(Q,\F)$ dans $\ha(\ha(\E,\opa),\F)$ est nulle.

\vs 0.4 cm

\th{L}{emme}{2.7} \textit{Pour tout $p$ et $q$ les diagrammes suivants
sont commutatifs :}
$$\begin{array}{ccc}
\ea^p(\ea^q(\E,\opa),\F) & \simeq & \ea^p(\ea^q(\E,\opa),\F)\\
\mapup{f_{p,q}} & & \mapup{g_p(\ea^q(\E,\opa))} \\
\ek^p(\ek^q(E,\opk),\F) & \fl &
\ek^p(\ea^q(\E,\opa)\ota\opk,\F)\\
\end{array}$$

\vs 0.2 cm

\dm :
Pour simplifier les notations, nous ne pr\'eciserons plus le faisceau
${\cal G}$ associ\'e au morphisme $g_p({\cal G})$. La fonctorialit\'e
de l'application de restriction donne le diagramme commutatif
suivant :
$$\begin{array}{ccc}
\ea^p(\ek^q(E,\opk),\F) & \stackrel{{\rm res}_A}{\fl} &
\ea^p(\ea^q(\E,\opa),\F)\\
\mapup{g_p} & & \mapup{g_p} \\
\ek^p(\ek^q(E,\opk),\F) & \stackrel{{\rm res}_k}{\fl} &
\ek^p(\ea^q(\E,\opa)\ota\opk,\F)\\
\end{array}$$
qui permet de conclure car $f_{p,q}$ est la compos\'ee de ${\rm
res}_A$ et de $g_p$.

\vs 0.2 cm

L'isomorphisme de $Q$ dans $\ea^1(\E,\opa)\ota\opk$ (car la limite est
r\'eduite) et le lemme pr\'ec\'edent nous donnent le diagramme
commutatif :
$$\begin{array}{ccc}
\ea^2(\ea^1(\E,\opa),\F) & \simeq & \ea^2(\ea^1(\E,\opa),\F)\\
\mapup{f_{2,1}} & & \mapup{g_2} \\
\ek^2(\ek^1(E,\opk),\F) & \stackrel{i}{\fl} &
\ek^2(Q,\F)\\
\end{array}$$

Nous savons que l'image du faisceau $\hk(Q,\F)$ dans
$\ea^2(\ea^1(\E,\opa),\F)$ est nulle. L'injectivit\'e de $g_2$ montre
que l'image de $\hk(Q,\F)$ dans $\ek^2(Q,\F)$ est nulle.

\subsection{Une nouvelle condition}

Dans ce paragraphe nous donnons une nouvelle condition limite.
Nous supposons que $E$ v\'erifie les
conclusions du th\'eor\`eme 2.6 et que son lieu singulier est de dimension
pure \'egale \`a 1.

Cette derni\`ere hypoth\`ese assure que le faisceau $E''$ est
localement libre et que le faisceau $\ek^1(E,\opk)$ s'identifie au
faisceau $Q$. Notons $C$ son support sch\'ematique.

\vs 0.4 cm

On consid\`ere la suite exacte $0\fl E\fl E''\fl Q\fl 0$ \`a laquelle
on applique le foncteur $\hk(\ \bullet\ ,\F)$. On obtient la
suite exacte longue suivante :
$$0\fl\hk(Q,\F)\stackrel{p}{\fl}\hk(E'',\F)\stackrel{q}{\fl}\hk(E,\F)
\fl\ek^1(Q,\F)\fl 0$$
Soit $F$ le conoyau de $p$.
La suite exacte longue pr\'ec\'edente d\'efini un \'el\'ement
$\xi'$ du groupe ${\rm Ext}_{\oo_C}^1(\ek^1(Q,\F),F)$.

Par ailleurs, nous savons que la fl\`eche de $\hk(Q,\F)$ dans
$\ek^2(Q,\F)$ (qui est la compos\'ee de $p$ et $f$) est nulle. Ceci
donne une fl\`eche de $F$ dans $\ek^2(Q,\F)$ qui permet
de d\'efinir, \`a partir de $\xi'$, un \'el\'ement $\xi\in{\rm
Ext}^1_{\oo_C}(\ek^1(Q,\F),\ek^2(Q,\F))$.

\vs 0.4 cm

\th{T}{h\'eor\`eme}{2.8} \textit{Si $E$ est limite r\'eduite de
fibr\'es vectoriels, alors on a l'annulation de l'\'el\'ement $\xi$ du
groupe ${\rm Ext}_{\oo_C}^1(\ek^1(Q,\F),\ek^2(Q,\F))$.}

\vs 0.2 cm

\dm :
On a le diagramme commutatif suivant (la seconde ligne est
scind\'ee, cf. lemme 1.10) :
$$\begin{array}{ccccccccc}
 & & \hk(E'',\F) & \to & \hk(E,\F) & \to & \ek^1(Q,\F) & \to & 0 \\
 & & \downarrow & & \downarrow & & \mapsdown{j} & & \\
0 & \to & \ek^2(Q,\F) & \to & \ea^2(\ea^1(\E,\opa),\F) & \to &
\ek^1(Q,\F) & \to & 0 \\
  \end{array}$$
Pour montrer le th\'eor\`eme, il suffit de montrer que $j$ est un
isomorphisme.

Etudions les noyaux et conoyaux des deux premi\`eres
fl\`eches verticales.
Sous nos hypoth\`eses, le complexe $\ECC$ n'a que deux termes et sa
cohomologie est $E''$ et $Q$
en degr\'es 1 et 0. On applique alors le foncteur $\hk(\ \bullet\
,\F)$ \`a ce complexe pour obtenir la suite spectrale :
$$\ ^2E^{p,q}=\ek^p(\ek^q(E,\opk),\F)\Rightarrow H^{p+q}(\hk(\ECC,\F))$$
Elle permet de montrer que le conoyau de la fl\`eche
$\hk(E'',\F)\fl\ek^2(Q,\F)$ est nul et que l'on a la suite exacte :
$$H^2(\hk(\ECC,\F))\fl\hk(E'',\F)\fl\ek^2(Q,\F)\fl 0$$

Le m\^eme raisonnement avec
le complexe $\ccE$ et le
foncteur $\ha(\ \bullet\ ,\F)$, donne :
$$H^2(\ha(\ccE,\F))\fl\ha(\E,\F)\fl\ea^2(\ea^1(\E,\opa),\F)\fl 0$$
Cependant les complexes $\hk(\ECC,\F)$ et $\ha(\ccE,\F)$ sont identiques
ainsi on compl\`ete le diagramme pr\'ec\'edent en :
$$\begin{array}{ccccccccc}
 & & H^2(\hk(\ECC,\F)) & \simeq & H^2(\ha(\ccE,\F)) &  & & & \\
 & & \downarrow & & \downarrow & & & & \\
 & & \hk(E'',\F) & \to & \hk(E,\F) & \to & \ek^1(Q,\F) & \to & 0 \\
 & & \downarrow & & \downarrow & & \mapsdown{j} & & \\
 & & \ek^2(Q,\F) & \hookrightarrow & \ea^2(\ea^1(\E,\opa),\F) & \to &
\ek^1(Q,\F) & \to & 0 \\
 & & \downarrow & & \downarrow & & & & \\
 & & 0 & & 0 & & & & \\
  \end{array}$$
ce qui montre que $j$ est un isomorphisme.

\vs 0.4 cm

\def \mod {{\mathfrak{Mod}}}
\def \ab {{\mathfrak{Ab}}}

Nous montrons maintenant qu'il suffit de v\'erifier les conditions des
th\'eor\`emes 2.6 et 2.8 pour le faisceau $\F=Q$.

Soit $B$ un anneau commutatif noeth\'erien et $\mod(B)$ la cat\'egorie des
$B$-modules de type fini. On dit qu'un foncteur additif $F:\mod(B)\fl\ab$
est coh\'erent si il admet une pr\'esentation repr\'esentable : $$h_N\fl
h_M\fl F\fl 0$$ o\`u $M$ et $N$ sont des objets de $\mod(B)$ et $h_N$ est
le foncteur ${\rm Hom}_{\mod(B)}(N,\bullet)$. Les foncteurs coh\'erents
forment une cat\'egorie ab\'elienne not\'ee $\coh(B)$. C'est une sous
cat\'egorie pleine de celle des foncteurs.

En prenant un recouvrement affine de $X$ par des ouverts $U_i$ d'anneau
$B_i$, nous pouvons d\'efinir la cat\'egorie $\coh(X)$ : les objets sont
les collections $(F_i,u_{i,j})$ o\`u $F_i$ est un \'el\'ement de
$\coh(B_i)$ et $u_{i,j}$ un isomorphime entre $F_i$ et $F_j$ sur l'ouvert
$U_{i,j}$. Les $u_{i,j}$ doivent \'evidement v\'erifier la condition de
cocycle sur $U_{i,j,k}$.

\vs 0.4 cm

\th{P}{roposition}{2.9} \textit{(\i) Sous les hypoth\`eses du th\'eor\`eme
2.6, on a l'identification :} $${\rm
Hom}_{\coh(X)}(\h(Q,\bullet),\e^2(Q,\bullet))\simeq\e^2(Q,Q)$$

\textit{(\i\i) Sous les hypoth\`eses du th\'eor\`eme 2.8, on a
l'identification :} $${\rm
Ext}^1_{\coh(X)}(\e^1(Q,\bullet),\e^2(Q,\bullet))\simeq H^1\e^1(Q,Q).$$

\vs 0.2 cm

\dm : Le premier cas vient directement du fait que l'on a (cf. \cite{Ha2}
lemme 1.2) $${\rm Hom}_{\coh(B)}(h_N,F)\simeq F(N).$$ Pour le second cas,
on commence par l'identification locale suivante :

\vs 0.2 cm

\th{L}{emme}{2.10} \textit{Pour tout $i$, on a} $${\rm
Hom}_{\coh(B_i)}(\ev^1_{B_i}(Q\vert_{U_i},\bullet),\ev^2_{B_i}(Q\vert_{U_i},\bullet))\simeq
\ev^1_{B_i}(Q\vert_{U_i},Q\vert_{U_i})$$

\vs 0.2 cm

\dm : Le faisceau $Q$ est ici de dimension projective 2 et il existe une
r\'esolution localement libre sym\'etrique de ce faisceau. Nous en
d\'eduisons que $\e^i(Q,\bullet)\simeq\t^{\oo_X}_{2-i}(Q,\bullet)$ pour
$i\in\{0,1,2\}$. Nous avons donc $${\rm
Hom}_{\coh(B_i)}(\ev^1_{B_i}(Q\vert_{U_i},\bullet),\ev^2_{B_i}(Q\vert_{U_i},\bullet))\simeq{\rm
Hom}_{\coh(B_i)}(\t_1^{B_i}(Q\vert_{U_i},\bullet),Q\vert_{U_i}\ot_{B_i}\bullet).$$

Par ailleurs, il existe une involution not\'ee $*$ de $\coh(B_i)$ dans lui
m\^eme qui est exacte contravariante et qui echange $\ev^i(M,\bullet)$ et
$\t^i(M,\bullet)$ (\cite{Ha2} paragraphe 4). Le groupe d'homomorphisme
pr\'ec\'edent est donc isomorphe \`a $${\rm
Hom}_{\coh(B_i)}(h_{Q\vert_{U_i}},\ev^1_{B_i}(Q\vert_{U_i},\bullet))
\simeq\ev^1_{B_i}(Q\vert_{U_i},Q\vert_{U_i}).$$

Pour $F$ et $G$ des objets de $\coh(X)$, nous pouvons d\'efinir un
$\oo_X$-module $\hv(F,G)$ qui est localement ${\rm
Hom}_{\coh(B_i)}(F_i,G_i)$. Le lemme 2.10 nous montre que
$$\hv(\e^1(Q,\bullet),\e^2(Q,\bullet))\simeq\e^1(Q,Q).$$ Pour conclure
nous utilisons la suite spectale de composition de foncteurs suivante :
$$H^p(\ev^q(F,G))\Rightarrow{\rm Ext}^{p+q}_{\coh(X)}(F,G).$$ On
l'applique aux foncteurs $F=\e^1(Q,\bullet)$ et $G=\e^2(Q,\bullet)\simeq
Q\ot\bullet$~. Le foncteur $G$ \'etant injectif (cf. \cite{Ha2} proposition
4.9), nous obtenons l'identification voulue. 

\vs 0.4 cm

 \thf{Remarques}{2.11} (\i) Le groupe ${\rm
Ext}^1_{\oo_C}(\e^1(Q,Q),\e^2(Q,Q))$ est exactement le groupe $H^1\e^1(Q,Q)$. Par
ailleurs, ce dernier groupe s'identifie \`a $H^1N_C$ o\`u $N_C$ est le
fibr\'e normal \`a $C$. Nous avons donc identifi\'e une obstruction dans
ce dernier groupe. C'est dans ce m\^eme groupe qu'appara\^\i t dans
\cite{ES} ce que nous appelons la condition d'Ellingsrud et Str\o mme.

(\i\i) Nous donnons ici un exemple lorsque $X$ n'est pas $\p^3$. Soit $X$
le volume de Fano cubique lisse de $\p^4$. Le faisceau inversible
$\omega_X$ est $\oo_X(-2)$ il admet donc $\oo_X(-1)$ comme racine
carr\'ee. L'espace des modules $\M_X(0,2,0)$ a \'et\'e \'etudi\'e dans
\cite{D} et \cite{IM}. En particulier, S. Druel donne une description
compl\`ete du bord de l'ouvert des fibr\'es vectoriels. Un faisceau $E$ de
ce bord est de l'une des deux formes suivantes : $$0\fl
E\fl\oo_X^2\fl\theta_C(1)\fl 0$$ o\`u $\theta_C$ est une th\'eta
caract\'eristique sur une conique lisse de $X$, ou bien : $$0\fl
E\fl\oo_X^2\fl\theta_{L_1}(1)\oplus\theta_{L_2}(1)\fl 0$$ o\`u
$\theta_{L_i}$ est une th\'eta carct\'eristique sur la droite $L_i$. Ces
faisceaux v\'erifient donc les conditions du th\'eor\`eme 2.3 mais aussi
celles des th\'eor\`eme 2.6 et 2.8 qui sont dans ce cas triviales.

On ne connait pas de d\'eformation explicite ayant pour limite les
faisceaux pr\'ec\'edents. On ne peut \`a fortiori en conna\^\i tre de
r\'eduite. Il semble cependant raisonnable de penser qu'il existe des
d\'eformations r\'eduites ayant pour limite les faisceaux pr\'ec\'edents
(c'est ce qu'affirme la conjecture 0.7). 

(\i\i\i) Nous donnons dans \cite{P2} un exemple complet lorsque $X$ est la
quadrique lisse de $\p^4$. Nous v\'erifions que les conditions sont
suffisantes et qu'une d\'eformation g\'en\'erale est r\'eduite.

\section{D\'eformations de faisceaux par homoth\'etie}

\def \ho {\underline{\bf Hom}}
\def \ex {\underline{\bf Ext}}
\def \to {\underline{\bf Tor}}
\def \ec {\widetilde{\p^3}}
\def \eca {\widetilde{\p^3_A}}
\def \eck {\widetilde{\p^3_k}}
\def \oo {{\cal O}}
\def \oec {{\cal O}_{\ec}}
\def \oeca {{\cal O}_{\widetilde{\p^3_A}}}
\def \oeck {{\cal O}_{\widetilde{\p^3_k}}}
\def \ecE {\widetilde{E}}
\def \ecaE {\widetilde{\E}}
\def \cM {{\cal M}}
\def \bV {{\bf V}}
\def \otau {\oo_{\ec}(\tau)}

Dans ce paragraphe, nous d\'ecrivons une famille $\bV$ de faisceaux de
$\M_{\p^3}(0,c_2,0)$ pour
laquelle les conditions pr\'ec\'edentes sont non triviales et suffisantes. Nous
construisons ainsi en particulier une famille adh\'erente \`a la
vari\'et\'e $\bI_{c_2}$ des instantons qui est toujours en codimension
1 et irr\'eductible pour $c_2\leq 5$.

Soit $\theta$ une th\'eta-caract\'eristique plane et soit $W$ sous-espace
vectoriel de dimension $2$ de sections de $\theta(2)$ tel que la
fl\`eche $W\ot\op\fl\theta(2)$ est surjective. Notons $\bV$ la famille
de $\M_{\p^3}(0,c_2,0)$ des faisceaux $E_0$ d\'efinis par :
\vs -0.3 cm
$$0\fl E_0\fl W\ot\op\fl\theta(2)\fl 0.$$
La dimension de $\bV$ est $\frac{1}{2}c_2(c_2+3)+4c_2-1$. La
dimension attendue (et obtenue au moins en bas degr\'es) de $\U$ est
$8c_2-3$.
On ne s'attend donc pas \`a ce que la vari\'et\'e $\bV$ soit adh\'erente
\`a $\U$.

Nous utilisons la construction de G. Ellingsrud et S.A. Str\o
mme \cite{ES} pour d\'ecrire des d\'eforma\-tions ``par homoth\'etie''
de fibr\'es vectoriels. Nous montrons
(th\'eor\`eme 3.2) qu'un ferm\'e $\bV_0$ de $\bV$ est adh\'erent \`a $\U$.
Il est d\'efini par les couples $(\theta,W)$ v\'erifiant ce que nous
appelons ``la condition d'Ellingsrud et St\o mme'', qui appara\^\i t
dans \cite{ES} et dont nous montrons qu'elle est impliqu\'ee par
la condition du th\'eor\`eme 2.8 (proposition 3.9).

Nous d\'ecrivons par ailleurs (th\'eor\`eme 3.2) un param\'etrage de
$\bV_0$  qui prolonge celui d'un ouvert $\U_P$ de $\U$ d\'ecrit dans
\cite{ES}. Nous v\'erifions ainsi que $\bV_0$ est de codimension 1 dans
$\U_P$.

\vs -0.6 cm

\subsection{Notations et rappels}

Nous rappelons ici quelques notations et r\'esultats de \cite{ES}.
Soit $P$ un point de $\p^3$ et soit
$p:\ec\fl\p^3$
l'\'eclatement de
ce point dans $\p^3$.
On a un morphisme $q:\ec\fl\p^2$ o\`u $\p^2$ est la vari\'et\'e des
droites de $\p^3$ passant par $P$.
Le morphisme $q$ est le fibr\'e en droites projectives associ\'e au
faisceau $\oo_{\p^2}\oplus\oo_{\p^2}(1)$. On note
$\oec(\sigma)$ (resp. $\oec(\tau)$) le faisceau $q^*\oo_{\p^2}(1)$
(resp. $p^*\oo_{\p^3}(1)$). Le diviseur exceptionnel est not\'e $B$.

\vs 0.2 cm

Soit $E$ un faisceau localement libre de rang 2 sur $\p^3$ et
supposons que
\begin{itemize}
\item il existe au moins une droite non sauteuse pour $E$ passant par
$P$,
\item il n'existe pas de droite bisauteuse pour $E$ passant par $P$.
\end{itemize}
Ces condition sont v\'erifi\'ees sur un ouvert non vide $\U_p$ de
$\U$.
On note $\ecE$ le faisceau $p^*E$, on a $E=p_*\ecE$.
On note $F$ le faisceau $q_*p^*E$ et $\theta$ le faisceau
$R^1q_*p^*E(-1)$. Le faisceau $\theta$ est une th\'eta
caract\'eristique, soit $C$ son support.
Le faisceau $F$ est d\'efini par deux
sections de $\theta(2)$ :
\vs -0.3 cm
$$0\fl F\fl\oo_{\p^2}^2\fl\theta(2)\fl 0\ \ \ (1)$$
Le faisceau $\ecE$ est le noyau
d'une fl\`eche surjective $q^* {\check F}
\stackrel{\phi}{\longrightarrow} q^* \theta (\sigma +\tau ).$

G.~Ellingsrud et S.A.~Str\o mme montrent que, pour tout faisceau $F$
d\'efini par deux sections de $\theta(2)$ o\`u $\theta$ est une
th\'eta-caract\'eristique plane, on a

\vs 0.2 cm

\th{P}{roposition}{\cite{ES}} \textit{Le faisceau $F$ admet la
d\'ecomposition ${\check F}\vert_C\simeq\theta(1)\oplus\theta(2)$
not\'ee $(D)$ si et seulement si il existe des fl\`eches surjectives
$q^* {\check F} \stackrel{\phi}{\longrightarrow} q^* \theta (\sigma +
\tau ).$}

\vs 0.2 cm

Ces conditions sont \'equivalentes \`a dire que $F$ provient
d'un faisceau $E\in\U$.
La d\'ecomposition
$(D)$ est ce qu'on appele la \textit{condition d'Ellingsrud et Str\o
mme}, elle a \'et\'e exprim\'ee directement sur la
th\'eta-caract\'eristique par \cite{HA} (voir plus loin), nous la
r\'einterpr\'etons \`a la proposition 3.9.

Soit $\theta$ une th\'eta-caract\'eristique plane et $F$ un faisceau
d\'efini par la suite exacte $(1)$. Il existe une fl\`eche
surjective $m_1$ unique \`a homoth\'etie pr\`es :
\vs -0.4 cm
$$0\fl\oo_{\p^2}^2\fl{\check F}\stackrel{m_1}{\fl}\theta(1)\fl
0\ \ \ (2).$$
\vs -0.1 cm
En particulier on a $\phi\vert_B=\lambda m_1$ o\`u $\lambda$ est un
scalaire non nul. La resctriction de $m_1$ \`a la courbe $C$ donne la
suite exacte :
\vs -0.4 cm
$$0\fl\theta(2)\fl{\check F}\vert_C\fl\theta(1)\fl 0\ \ \ (3).$$
\vs -0.1 cm

\subsection{Les familles $\bV$ et $\bV_0$}

G. Ellingsrud et S.A Str\o mme ont donn\'e dans \cite{ES} un
param\'etrage birationnel de l'ouvert $\U_P$. Rappelons le rapidement.

Notons
${\rm {\bf Theta}}_{\p^2}(c_2)$ le sch\'ema des
th\'eta-caract\'eristiques planes sur une courbe de degr\'e $c_2$
(suppos\'ees localement libre sur leur support). Sur ${\rm {\bf
Theta}}_{\p^2}(c_2)\times\p^2$, il existe une th\'eta-caract\'eristique
universelle $\Theta$.  Consid\'erons la grassmannienne ${\rm
Grass}(2,{pr_1}_*(\Theta(2))^\vee)$ des sous-fibr\'es de rang 2 de
$pr_1(\Theta(2))^\vee$. On note $K$ le sous fibr\'e
tautologique et ${\bf G}$ l'ouvert de la grassmannienne o\`u le
morphisme naturel
$$K\boxtimes\oo_{\p^2}\fl\Theta(2)$$
est surjectif. On note alors $F$ le noyau de cette fl\`eche.

Pour qu'un tel faisceau provienne d'un fibr\'e sur $\p^3$,
G. Ellingsrud et S.A. Str\o mme ont identifi\'e une condition de
scindage $(D)$ (cf. proposition [ES]). F. Han \cite{HA} a \'etudi\'e
cette condition qu'il traduit de la mani\`ere suivante : pour toute
th\'eta-caract\'eristique plane $\theta$, il existe une fl\`eche
$$\Lambda^2H^0\theta(2)\stackrel{\psi}{\fl}H^1\oo_C(1).$$
Le faisceau $F$ est d\'efini par deux sections de $\theta(2)$
c'est \`a dire par un point $\xi'$ de $\Lambda^2H^0\theta(2)$. F.~Han
montre alors que $F$ provient d'un fibr\'e vectoriel sur $\p^3$ si
et seulement si $\psi(\xi')$ est nul.

Notons ${\bf G_0}$ le ferm\'e de ${\bf G}$ d\'efini par cette
derni\`ere condition. Sur ${\bf G_0}$ on d\'efinit alors le faisceau
$R={pr_{\bf G_0}}_*(\ho_{\ec\times{\bf G_0}}(q^*{\check
F},q^*\Theta(\sigma+\tau)))$.
Le noyau de la fl\`eche universelle de $q^*{\check F}$ dans
$q^*\Theta(\sigma+\tau)$ d\'efini une application rationelle
$\pi:\p_{\bf G_0}({\check R})\dasharrow\M_{\p^3}(0,n,0)$.

\vs 0.4 cm

\th{T}{h\'eor\`eme}{[ES]} \textit{L'image par $\pi$ de l'ouvert de $\p_{\bf
G_0}({\check R})$ correspondant aux fl\`eches de $q^*{\check F}$ dans
$q^*\Theta(\sigma+\tau)$ surjectives est isomorphe \`a
$\U_P$.}

\vs 0.4 cm

\th{D}{\'efinition}{3.1} Nous d\'efinissons la vari\'et\'e $\bV$
(resp. $\bV_0$) en faisant la construction de ${\bf G}$ (resp. ${\bf
G}_0$) en famille au dessus de la vari\'et\'e des plans de $\p^3$. Le
noyau de la fl\`eche universelle surjective $\oo_{\bV}^2\fl\Theta(2)$
d\'efini une immersion de $\bV$ (et donc de $\bV_0$) dans
$\M_{\p^3}(0,c_2,0)$.

\vs 0.4 cm

La suite exacte de d\'efinition de $F$ donne par dualit\'e une
fl\`eche
$$m_1:q^*{\check F}\fl q^*\Theta(\sigma)$$
qui est unique \`a homoth\'etie pr\`es. Ainsi le faisceau ${pr_{{\bf
G}_0}}_*(\ho_{\ec\times{\bf G}_0}(q^*{\check F},q^*\Theta(\sigma)))$ est
inversible, notons le $\L$. Le faisceau ${pr_{{\bf
G}_0}}_*(\ho_{\ec\times{\bf G}_0}(\ec,\otau))$ est isomorphe au
faisceau trivial $H^0\op(1)\ot\oo_{{\bf G}_0}$.

Notons $R_0$ le faisceau $H^0\op(1)\ot\L$. En prenant l'image directe
par $pr_{{\bf G}_0}$ de la compos\'ee :
$$\ho_{\ec\times{\bf
G}_0}(\ec,\otau)\ot\ho_{\ec\times{\bf
G}_0}(q^*{\check
F},q^*\Theta(\sigma))\fl\ho_{\ec\times{\bf
G}_0}(q^*{\check
F},q^*\Theta(\sigma+\tau))$$
nous obtenons un fl\`eche injective de $R_0$ dans $R$.

\vs 0.4 cm

\th{T}{h\'eor\`eme}{3.2} \textit{La vari\'et\'e $\bV$ est de dimension
$\frac{1}{2}c_2(c_2+3)+4c_2-1$ et est contenue dans
$\M_{\p^3}(0,c_2,0)$. La vari\'et\'e
$\bV_0$
s'identifie \`a l'image par $\pi$ de $\p_{{\bf G}_0}({\check R_0})$, elle
est adh\'erente \`a $\U_P$ et de codimension 1.}

\textit{\`A un point $(\theta,W,H)$ de $\bV$ correpond un faisceau
$E_0$ :
$$0\fl E_0\fl W\ot\op\fl\theta_H(2)\fl 0$$
o\`u $\theta$ est une th\'eta caract\'eristique sur le plan $H$ et
$W\subset H^0\theta_H(2)$ est de dimension 2.}

\vs 0.2 cm

\dm :
Si $(\theta,W,H)\in\bV$,
le faisceau $E_0$ correpondant
est \'evidement d\'efini par la suite exacte ci-dessus. On calcule
facilement la dimension de $\bV$. Nous d\'ecrivons au prochain paragraphe
l'identification entre $\bV_0$ et $\pi(\p_{{\bf G}_0}({\check R_0}))$
ce qui prouve que $\bV_0$ est en codimension 1.

\vs 0.4 cm

\thf{Remarques}{3.3} (\i) Pour $c_2\leq 11$, F. Han
\cite{HA} remarque en citant \cite{ES} paragraphe (4.2) que la
vari\'et\'e ${\bf G}_0$ domine la vari\'et\'e ${\rm {\bf
Theta}}_{\p^2}(c_2)$. Cette derni\`ere a deux
(resp. trois) composantes irr\'eductibles si $c_2$ est pair
(resp. impair) \cite{S}. Ces composantes correspondent aux
th\'eta-caract\'eristiques paires impaires et canoniques (quand elles
existent). Nous notons $\bV_0^{pair}$, $\bV_0^{impair}$ et $\bV_0^{can}$
les composantes de $\bV_0$ associ\'ees.

Soit $E_0$ un faisceau de $\bV_0$ et $\theta$ la
th\'eta-caract\'eristique correspondante. Si $E\in\U_P$ se d\'eforme
en
$E_0$, on a l'\'egalit\'e
$h^0\theta\equiv h^1E(-2)$ (mod 2). En d'autres termes, la parit\'e de
$\theta$ est \'egale \`a celle de l'invariant d'Atiyah-Rees de $E$.

(\i\i) Si $c_2\leq 5$, F. Han a montr\'e que ${\bf G}_0^{pair}$ est
irr\'eductible. La famille $\bV_0^{pair}$ forme donc une composante
irr\'eductible de $\bI_{c_2}$.

\subsection{La d\'eformation}

Notons $A$ l'anneau de valuation discr\`ete $k[a]_{(a)}$. Pour tout
faisceau $M$ d\'efini sur ${\rm Spec}(k)$, nous noterons $M_A$ le
faisceau $M\ot_kA$. Nous nous donnons un fibr\'e
$E\in\M_{\p^3}(0,c_2,0)$ et un point $P$ v\'erifiant les conditions du
paragraphe 3.1. Nous construisons une d\'eformation de $E$
vers un point de $\bV_0$.

Soit $\H$ le transform\'e strict d'un hyperplan $H$ de $\p^3$ ne
rencontrant pas $P$. Notons encore $H\in H^0\oec(\tau)$ l'\'equation de $\H$
dans $\ec$.
Notons $\phi_A$ la fl\`eche $a\phi+H\circ q^*m_1$ (c'est une
d\'eformation d'un point g\'en\'eral de $\p_{\bf G_0}({\check R})$
vers un point de $\p_{\bf G_0}({\check R}_0)$). Nous construisons la
famille $\ecaE$ suivante :
$$0\fl\ecaE\fl q^*{\check F}_A \stackrel{\phi_A}{\fl}
q^*\theta_A(\sigma+\tau)$$
La platitude de $q^*{\check F}_A$ impose celle de $\ecaE$. En
appliquant le foncteur $p_*$, on a la suite exacte
$$0\fl p_*\ecaE\fl p_*q^*{\check F}_A \stackrel{\phi_A}{\fl}
p_*q^*\theta_A(\sigma+\tau)$$
qui montre que la famille $\E=p_*\ecaE$ est plate. Le faisceau
g\'en\'eral de $\E$ est un instanton. Le morphisme $q$
identifie $H$ \`a $\p^2$. Nous noterons
$\theta_H(2)$ le faisceau
$q^*\theta(\sigma+\tau)\vert_H\simeq\theta(2)$. Pour terminer la
d\'emonstration du th\'eor\`eme 3.2, nous devons montrer que la limite
$E_0$ est un \'el\'ement de $\bV_0$. Nous d\'ecrivons ainsi
l'identification entre $\bV_0$ et $\p_{\bf G_0}({\check R}_0)$.

\vs 0.4 cm

\th{P}{roposition}{3.4} \textit{Le conoyau de $\phi_A$ est $\theta_H(2)$.}

\vs 0.2 cm

\dm :

\th{L}{emme}{3.5} \textit{La fl\`eche $f_A$ suivante est surjective :}
\vs -0.2 cm
$$\ q^*{\check F}_A\oplus
q^*\theta_A(\sigma)\stackrel{f_A}{\longrightarrow}q^*\theta_A(\sigma
+\tau)\oplus q^*\theta_A(\sigma)\ \ \ {\rm avec}\ \ \
f_A={\begin{pmatrix}
\phi & H \\
m_1 & -a
\end{pmatrix}}.$$

\dm :
Notons $f_0$ la sp\'ecialisation de $f_A$ sur le corps r\'esiduel. Il
suffit de montrer la surjectivit\'e de $f_0$ (lemme de Nakayama). Elle
est \'equivalente \`a la surjectivit\'e de la fl\`eche :
$$q^*{\check F} \stackrel{\binom{\phi\vert_H}{m_1}}{\longrightarrow}
\theta_H(2)\oplus q^*\theta(\sigma)$$
Cette surjectivit\'e ne pose probl\`eme que sur le support des
faisceaux \`a l'arriv\'ee. De plus, comme $\phi$ et $m_1$ sont
surjectives, la question de la surjectivit\'e se pose seulement sur
$H\cap S\simeq C$. On a :

$$\theta(1)\oplus\theta(2)\simeq{\check
F}\vert_{C}\stackrel{\binom{\phi\vert_H}{m_1}}{\longrightarrow}
\theta(2)\oplus\theta(1)$$
Les fl\`eches $m_1$ et $\phi\vert_H$ prises ind\'ependement sont
surjectives. Ainsi les deux endomorphismes de $\theta(1)$ et
$\theta(2)$ sont des scalaires non nuls et on a la
surjectivit\'e. 

\vs 0.2 cm

La fl\`eche $\phi_A$ s'inscrit dans le diagramme suivant :
$$\begin{array}{ccc}
 q^*{\check F}_A\oplus  q^*\theta_A(\sigma) &
\stackrel{f_A}{\longrightarrow} & q^*\theta_A(\sigma +\tau)\oplus
q^*\theta_A(\sigma) \\
\downarrow &   & \mapdown{(a,H)} \\
q^*{\check F}_A & \stackrel{\phi_A}{\longrightarrow} &
q^*\theta_A(\sigma+\tau) \\
\end{array}$$

ce qui nous permet de dire gr\^ace au lemme 3.5 que le conoyau de
$\phi_A$ est celui de la fl\`eche
$$q^*\theta_A(\sigma+\tau)\oplus
q^*\theta_A(\sigma)\stackrel{(a,H)}{\fl}q^*\theta_A(\sigma+\tau).$$

Par ailleurs, la suite exacte suivante tensoris\'ee par
$q^*\theta_A(\sigma)$ reste exacte :
$$0\longrightarrow \oeca \stackrel{\binom{H}{-a}}{\longrightarrow}
\oeca(\tau)\oplus\oeca \stackrel{(a,H)}{\longrightarrow}
\oeca(\tau)\longrightarrow \oo_H(1)\ot_{\oeca}\oeck\fl 0\ \ \ (4).$$

\vs 0.4 cm

\th{P}{roposition}{3.6} \textit{Le faisceau $R^1p_*\ecaE$ est nul.}

\vs 0.2 cm

\dm :
Calculons la restriction au diviseur exceptionnel $B$ de la suite
exacte de d\'efinition de $\ecaE$. Rappelons que $\phi\vert_B=\lambda
m_1$. Nous savons que $H\cap B$ est vide et que l'\'equation $H$ est
inversible sur $B$. On a la suite exacte :
$$0\fl\ecaE\vert_B\fl{\check
F}_A\stackrel{(a\lambda+H)m_1}{\fl}\theta_A(1)\fl 0.$$
Le scalaire $a\lambda+H$ est inversible donc le faisceau
$\ecaE\vert_B$ est isomorphe \`a $\oo_{\p^2_A}^2$ (cf. suite exacte
$(2)$). On montre alors que le faisceau $R^1p_*\ecaE$ est nul (voir
par exemple \cite{ES} proposition (2.4)). 

\vs 0.4 cm

\th{T}{h\'eor\`eme}{3.7} \textit{Le faisceau limite $E_0$ de la famille
$\E$ est un faisceau sans torsion non localement libre. Son bidual
$E_0''$ est trivial et le conoyau de l'injection canonique de $E_0$
dans $E_0''$ est $\theta(2)$ o\`u $\theta$ est une
th\'eta-caract\'eristique sur une courbe plane.}

\vs 0.2 cm

\dm:
La proposition 3.6 montre que, si on note $\ecE_0$
(resp. $(p_*\ecaE)_0$) la sp\'ecialisation au corps r\'esiduel du
faisceau $\ecaE$ (resp. $(p_*\ecaE)$), alors on a
$(p_*\ecaE)_0=p_*\ecE_0$. C'est le faisceau limite, notons le $E_0$.

Par ailleurs, la suite exacte $(4)$ montre que la restriction de
l'image de $\phi_A$ au corps r\'esiduel est
$q^*\theta(\sigma)\oplus\theta_H(2)$. La fl\`eche de ${\check F}$ dans
$q^*\theta(\sigma)$ est $m_1$

ce qui nous donne la suite exacte
$$0\longrightarrow E_0 \longrightarrow \op^2 \longrightarrow
\theta_H(2)\longrightarrow 0\ \ \  (5).$$
Ceci termine \'egalement la d\'emonstration du th\'eor\`eme 3.2.

\vs 0.4 cm

\thf{Remarque}{3.8} La d\'eformation $\E$ est r\'eduite. En effet, le
faisceau $p^*\ea^1(\E,\opa)$ est $\ex_{\oeca}^1(\ecaE,\oeca)$. La suite
exacte de d\'efinition de $\ecaE$ montre que
$\ex_{\oeca}^1(\ecaE,\oeca)$ s'identifie \`a
$\ex_{\oeca}^3(\theta_H(2),\oeca)$. Il est donc annul\'e par $a$.

\vs 0.4 cm

Nous montrons maintenant que la \textit{condition d'Ellingsrud et
Str\o mme} $(D)$ est impliqu\'ee par le th\'eor\`eme 2.8.
On applique le
foncteur $\ho_{\op}(\ \bullet\ ,\theta(2))$ \`a la suite exacte
$(5)$. On obtient la suite exacte longue :
$$0\fl\ho_{\op}(\theta(2),\theta(2))\fl\theta(2)^2\fl
\ho_{\op}(E_0,\theta(2))\fl\ex_{\op}^1(\theta(2),\theta(2))\fl 0$$
qui d\'efini une extension $\xi$ de
$\ex_{\op}^1(\theta(2),\theta(2))$ par $\omega_C(4)$. Remarquons
que le faisceau $\omega_C(4)$ est isomorphe au faisceau
$\ex_{\op}^2(\theta(2),\theta(2))$.

\vs 0.4 cm

\th{P}{roposition}{3.9} \textit{La condition $(D)$ est \'equivalente \`a
l'annulation de l'\'el\'ement}
$$\xi\in{\rm Ext}^1(\ex_{\op}^1(\theta(2),\theta(2)),
\ex_{\op}^2(\theta(2),\theta(2))).$$

\vs 0.2 cm

\dm :
On identifie $H$ \`a $\p^2$. On applique le foncteur $\ho_{\oo_H}(\
\bullet\ ,\theta(2))$ \`a la suite exacte $(1)$. On obtient la
suite exacte longue :
$$0\fl\ho_{\oo_H}(\theta(2),\theta(2))\fl\theta(2)^2\fl
\ho_{\op}(F,\theta(2))\fl\ex_{\oo_H}^1(\theta(2),\theta(2))\fl 0$$
qui d\'efini une extension $\xi'$ de
$\ex_{\oo_H}^1(\theta(2),\theta(2))$ par $\omega_C(4)$. La condition
$(D)$ dit exactement que cette extention $\xi'$ est
triviale. Cependant on a l'identification
$$\ex_{\op}^1(\theta(2),\theta(2))\simeq
\ex_{\oo_H}^1(\theta(2),\theta(2))\oplus
\ho_{\oo_H}(\to^{\op}_1(\theta(2),\oo_H),\theta(2)).$$
Le faisceau $\ho_{\oo_H}(\to^{\op}_1(\theta(2),\oo_H),\theta(2))$ est
$\omega_C(3)$. Comme le groupe $H^1\omega_C(3)$ est nul, les
extensions $\xi$ et $\xi'$ sont nulles en m\^eme temps. Cette proposition
et le th\'eor\`eme 2.8 montrent en particulier que, dans $\bV$, seuls les
points de $\bV_0$ sont adh\'erents \`a $\U$. 

\vs 0.4 cm

\thf{Remarque}{3.10} En rempla\c cant le plan $H$ par le diviseur
exceptionnel $B$ et la fl\`eche $m_1$ par une fl\`eche $m_2$ (non
unique) de $q^*{\check F}$ dans $q^*\theta(2\sigma)$, on peut
construire de la m\^eme fa\c con des faisceaux $E_0$ limites de
fibr\'es vectoriels de la forme :
$$0\fl E_0\fl E''_0\fl Q\fl 0$$
o\`u $Q$ est un faisceau support\'e au point $P$ et $E''_0$ est
isomorphe \`a $p_*q^*N$ pour un faisceau $N$ localement libre sur
$\p^2$.

\vs 0.4 cm

\thf{Remarque}{3.11} Soit $\theta$ une th\'eta-caract\'eristique telle que
$\theta(2)$ est engendr\'e par ses sections. De la m\^eme fa\c con que
\cite{HA} nous pouvons construire une fl\`eche
$$\Lambda^2H^0\theta(2)\fl{\rm
  Ext}^1(\ex^1_{\op}(\theta(2),\theta(2)),\ex^2_{\op}(\theta(2),\theta(2)))$$
qui correspond \`a la condition du th\'eor\`eme 2.8 lorsque $E''$ est
trivial.

Remarquons que l'on a les identifications
$\ex^2_{\op}(\theta(2),\theta(2))\simeq\ho_{\oo_C}(\oo_C,\omega_C\ot\omega_X)$ et

\noi
$\ex^1_{\op}(\theta(2),\theta(2))\simeq\ho_{\oo_C}(\t^{\op}_1(\theta(2),\theta(2)),\omega_C\ot\omega_X)$, le groupe
d'extensions pr\'ec\'edent est donc isomorphe \`a
$H^1(\to^{\op}_1(\theta(2),\theta(2)))$.
Nous alors d\'eterminer une fl\`eche de $\Lambda^2H^0\theta(2)$ dans
$H^1\to^{\op}_1(\theta(2),\theta(2))$. Soit $L_0$ le faisceau
$H^0\theta(2)\ot\op$, on a le diagramme :
$$\begin{array}{ccccccccc}
0 & \fl & \alpha & \fl & L_0 & \fl & \theta(2) & \fl & 0 \\
 & & \uparrow & & \uparrow & & \Vert & & \\
0 & \fl & E & \fl & \op^2 & \fl & \theta(2) & \fl & 0
\end{array}$$
En tensorisant par $\theta(2)$ on a :
$$\begin{array}{ccccccccccc}
0 & \fl & \to^{\op}_1(\theta(2),\theta(2)) &
\fl & \alpha\ot\theta(2) & \fl & L_0\ot\theta(2) & \stackrel{f}{\fl} & \theta(2)\ot\theta(2)
& \fl & 0 \\
 & & \uparrow & & \uparrow & & \uparrow & & \Vert & & \\
0 & \fl & \to^{\op}_1(\theta(2),\theta(2)) &
\fl & E\ot\theta(2) & \fl & \op^2\ot\theta(2) & \stackrel{g}{\fl} & \theta(2)\ot\theta(2) &
\fl & 0
\end{array}$$
Notons $K$ et $K'$ les noyaux de $f$ et $g$. On a $K'\simeq\oo_C$ et le
diagramme :
$$\begin{array}{ccc}
H^0K & \fl & H^1\to^{\op}_1(\theta(2),\theta(2)) \\
\uparrow & & \Vert \\
H^0K' & \fl & H^1\to^{\op}_1(\theta(2),\theta(2))
\end{array}$$
Nous montrons que l'on peut factoriser la fl\`eche verticale de gauche par
$H^0\Lambda^2L_0$. On regarde le diagramme de complexes horizontaux :
$$\begin{array}{ccccccccc}
0 & \fl & \Lambda^2L_0 & \fl & L_0\ot\theta(2) & \fl & S^2\theta(2) & \fl &
0 \\
 & & \uparrow & & \uparrow & & \Vert & & \\
0 & \fl & \Lambda^2\op^2 & \fl & \op^2\ot\theta(2) & \fl & S^2\theta(2) & \fl
& 0
\end{array}$$
Ces deux complexes sont exacts sauf \`a gauche. En effet ce sont des
facteurs directs des complexes
$$0\fl L_0\ot L_0\fl
(L_0\ot\theta(2))^2\fl\theta(2)\ot\theta(2)\fl 0$$
$${\rm et}\ \ \ \ 0\fl\op\ot\op\fl
(\op^2\ot\theta(2))^2\fl\theta(2)\ot\theta(2)\fl 0$$
dont la cohomologie est
nulle sauf \`a gauche. Nous avons alors le diagramme :
$$\begin{array}{ccccccccc}
0 & \fl & {\rm Ker}_1 & \fl & \Lambda^2L_0 & \fl & K & \fl &
0 \\
 & & \uparrow & & \uparrow & & \Vert & & \\
0 & \fl & {\rm Ker}_2 & \fl & \Lambda^2\op^2 & \fl & K' & \fl &
0 \\
\end{array}$$
ce qui nous donne la factorisation recherch\'ee. La fl\`eche
$$\Lambda^2H^0\theta(2)=H^0\Lambda^2L_0\fl H^0K\fl
H^1\to^{\op}_1(\theta(2),\theta(2))$$
d\'efini la condition du th\'eor\`eme 2.8.

\vs 0.2 cm

\noi
\textsc{Mathematisches Institut der Universit\"at zu K\"oln}
\vs -0.1 cm
\noi
Weyertal 86-90
\vs -0.1 cm
\noi
D-50931 K\"oln
\vs -0.1 cm
\noi
email : \texttt{nperrin@mi.uni-koeln.de}


\begin{thebibliography}{99}


\bibitem[D]{D} \textit{St\'ephane Druel} : Espace des modules des
faisceaux de rang 2 semi-stables de classes de Chern $c\sb 1=0, c\sb
2=2$ et $c\sb 3=0$ sur la cubique de $\p\sp 4$,
Internat. Math. Res. Notices (2000) no. 19.

\bibitem[E]{Ei} \textit{David Eisenbud} : Commutative algebra. With a
view toward algebraic geometry, Graduate Texts in Mathematics,
150. Springer-Verlag, New York (1995).

\bibitem[ES]{ES} \textit{Geir Ellingsrud et Stein Arild Str\o mme} :
Stable rank-$2$ vector bundles on $\p^3$ with $c_1=0$ and $c_2=3$,
Math. Ann. 255 (1981).

\bibitem[H]{HA} \textit{Fr\'ed\'eric Han} : Sch\'ema des multisauteuses
d'un 4- ou 5-instanton et espace des modules. Th\`ese \`a
l'universit\'e de Lille (1996).

\bibitem[Ha1]{Ha} \textit{Robin Hartshorne} : Stable reflexive sheaves,
Math. Ann. 254 (1980).

\bibitem[Ha2]{Ha2} \textit{Robin Hartshorne} : Coherent functors,
  Adv. Math. 140 (1998), no. 1.  

\bibitem[IM]{IM} \textit{Atanas Iliev et Dimitri Markushevich} : The
Abel-Jacobi map for a cubic threefold and periods of Fano threefolds
of degree 14. Doc. Math. 5 (2000).

\bibitem[M]{MA} \textit{Masaki Maruyama} : Moduli of stable
sheaves. II, J. Math. Kyoto Univ. 18 (1978), no. 3.

\bibitem[NT]{NT} \textit{Mudumbai S. Narasimhan et Gunther Trautmann} :
Compactification of $M\sb {P\sb 3}(0,2)$ and Poncelet pairs of conics,
Pacific J. Math. 145 (1990), no. 2.

\bibitem[P1]{PE1} \textit{Nicolas Perrin} : Deux composantes du bord
des instantons de degr\'e 3, pr\'epublication disponible sur
math.AG/9901011 v2 (2001).

\bibitem[P2]{P2} \textit{Nicolas Perrin} : Limites de fibr\'es vectoriels
  dans $\M_{\Q_3}(0,2,0)$, en pr\'eparation.

\bibitem[S]{S} \textit{Christoph Sorger} : Th\^eta-caract\'eristiques
des courbes trac\'ees sur une surface lisse. J. Reine Angew. Math. 435
(1993).

%

\end{thebibliography}
\end{document}